\documentclass[11pt]{article}
\usepackage[dvips]{graphics,graphicx}
\DeclareGraphicsExtensions{.eps}

\newcounter{ctheorem}
\newcounter{clemma}
\newcounter{cdefinition}
\newtheorem{theorem}[ctheorem]{Theorem}
\newtheorem{lemma}[clemma]{Lemma}

\newtheorem{definition}[cdefinition]{Definition}

\def\ds{\displaystyle}
\def\id{\mbox{\bf 1}}
\def\hf{\left(\frac{+1}{2!}\right)}
\def\nn{\nonumber}
\def\re{\mbox{\sf Re}}

\def\stack{\stackrel}
\def\bspace{\!\!\!\!\!\!\!\!\!\!}
\def\bbspace{\!\!\!\!\!\!\!\!\!\!\!\!\!\!\!\!\!\!\!\!}
\def\bbbspace{\!\!\!\!\!\!\!\!\!\!\!\!\!\!\!\!\!\!\!\!\!\!\!\!\!\!\!\!\!\!}
\def\R{\bf\sf R}
\def\Z{\bf\sf Z}
\def\C{\bf\sf C}
\def\qed{\hbox{${\vcenter{\vbox{                        %HOLLOW SQUARE
   \hrule height 0.4pt\hbox{\vrule width 0.4pt height 6pt
   \kern5pt\vrule width 0.4pt}\hrule height 0.4pt}}}$}}

\begin{document}

\begin{center}
{\bf\Large Generalization of a relation between the Riemann zeta
  function and Bernoulli numbers}

\medskip

{\large S.C. Woon}

Trinity College, University of Cambridge, Cambridge CB2 1TQ, UK

s.c.woon@damtp.cam.ac.uk

MSC-class Primary 11M06; Secondary 11B68


December 23, 1998

\end{center}

\begin{abstract}
A generalization of a well-known relation between the Riemann zeta
function and Bernoulli numbers is obtained.  The formula is a new
representation of the Riemann zeta function in terms of a nested
series of Bernoulli numbers.
\end{abstract}

\section{A New Representation of the Riemann Zeta Function}

\begin{theorem}
\label{t:zetarep}
\begin{equation} \fbox{$
\begin{array}{lll} \zeta(s)&\!\!\!=\!\!\!&\ds -\; \frac{(2\pi)^s}{2} 
  \,w^{s-1}\!\!
  \lim_{\;{\ds \hat{s}}\to {\ds s}} \left\{\!
  \frac{ \left( \ds \frac{1}{2} + \sum_{n=1}^\infty (-1)^n
  {\hat{s}\!-\!1 \choose n}\!\!\left[ \frac{1}{2} +
  \sum_{m=1}^n \left(\frac{-1}{w}\right)^{\!\!m} \!\!  {n \choose m}
  \frac{B_{m+1}}{(m\!+\!1)!} \right] \!\right)}
{\ds \cos\left(\frac{\pi \hat{s}}{2}\right)} \right\}\end{array} $}
\label{e:zeta(s)Bn}
\end{equation}
for $\re(s)>(1/w)$ where $s\in\C, \;w\in\R, \;w>0$,
the notation of binomial coefficient is extended such that
$${s\!-\!1 \choose n} = \frac{1}{n!} \left[ \prod_{k=0}^{n-1} 
(s\!-\!1\!-\!k) \right] = \frac{1}{n!}\,
\frac{\Gamma(s)}{\Gamma(\!s\!-\!n)}\,,$$
$B_m$ are the Bernoulli numbers with $B_1 = 1/2$, and the limit only
needs to be taken when $s\in\{1,3,5,\dots\}$ for which the denominator
$\cos\,(\pi s/2)$ is $0$.
\end{theorem}

The representation (\ref{e:zeta(s)Bn}) can be seen as a generalization
of the well-known relation
$$\zeta(2n) \;=\;-\;\frac{(2\pi)^{2n}}{2} \,\frac{(-1)^n \,B_{2n}}{(2n)!}
\quad (n\in\Z^+)\;.$$

This representation (\ref{e:zeta(s)Bn}) of $\zeta(s)$ in terms of a
nested series of $B_n$ is distinct from the well-known Euler-Maclaurin
summation representation \cite[p.807, (23.2.3)]{abramowitz}
which also relates $\zeta(s)$ to $B_n$ as follows:
\begin{eqnarray} \zeta(s)&=&\lim_{N\to\infty}
\left[ \begin{array}{l}
\ds \sum_{n=1}^N n^{-s} - \frac{1}{-s\!+\!1}\,N^{-s+1} 
- \frac{1}{2}\,N^{-s}\\
\ds\Bigg. -\, \sum_{k=1}^M \frac{B_{2k}}{(2k)!}
\frac{\partial^{2k-1}}{\partial x^{2k-1}} x^{-s}\Big|_{x=N}\\
\ds\Big. +\, O(N^{-s-2M-1})
\end{array} \right]
\label{e:EulerMaclaurinsum}\\
&&\ds\bigg.(\re(s)>-2M-1)\;.\nn
\end{eqnarray}

To prove Theorem \ref{t:zetarep}, we shall have to introduce a binary
tree and a set of operators for generating Bernoulli Numbers.

\section{A Binary Tree for Generating Bernoulli Numbers}

\begin{definition}
\label{d:Bn}
Bernoulli numbers $B_n$ are defined by \cite[p.35, (1.13.1)]{bateman} 
\begin{equation}
\frac{z}{e^z-1} \,=\, \sum_{n=0}^\infty \frac{B_n}{n!} z^n \quad
(|z| < 2\pi)\;.\label{e:defBn}
\end{equation}
\end{definition}
Expanding the left hand side as a series and matching the coefficients
on both sides give
\begin{equation}
B_1 = -1/2, \quad B_n \left\{ \begin{array}{ccc} = 0 &,& 
\mbox{odd }n, \;n\ne 1\\ \ne  0 &,& \mbox{even }n \end{array} \right.\;.
\end{equation}
Now (\ref{e:defBn}) can be rewritten as
\begin{equation}
\frac{z}{e^z-1} + \frac{z}{2} \,=\, \sum_{n=0}^\infty
\frac{B_{2n}}{(2n)!} z^{2n}\;.\label{e:defB2n}
\end{equation}
Alternatively, $B_n$ can be defined as the solution of the recurrence
relation
\begin{equation}
 B_n \,=\, -\, \frac{1}{n+1} \,\sum_{k=0}^{n-1}
{n\!+\!1 \choose k} B_k\;, \quad B_0 = 1 \;.\label{e:defBnrecursion}
\end{equation}

A binary tree for generating Bernoulli numbers $B_n$ can be
constructed using two operators, $O_L$ and $O_R$.

\begin{definition}[{\bf The $B_n$-Generating Tree}]\quad\\
  At each node of the binary tree sits a formal expression of the form
  $\ds \frac{\pm 1}{a!\,b!\dots}$. The operators $O_L$ and $O_R$ are
  defined to act only on formal expressions of this form at the nodes
  of the tree as follows:
\end{definition}
\begin{eqnarray*} O_L^{} & : & \frac{\pm 1}{a!\,b!\dots} \to 
 \frac{\pm 1}{(a+1)!\,b!\dots} \;,\\
 O_R^{} & : & \frac{\pm 1}{a!\,b!\dots} \to 
 \frac{\mp 1}{2!\,a!\,b!\dots} \;.\end{eqnarray*}
Schematically, \begin{itemize}
\item $O_L^{}$ acting on a node of the tree generates a branch downwards to
the left (hence the subscript $L$ in $O_L^{}$) with a new node at the
end of the branch.
\item $O_R^{}$ acting on the same node generates a branch downwards to the
right. \end{itemize}

\begin{figure}[hbt]
\begin{center}
\includegraphics[height=175pt]{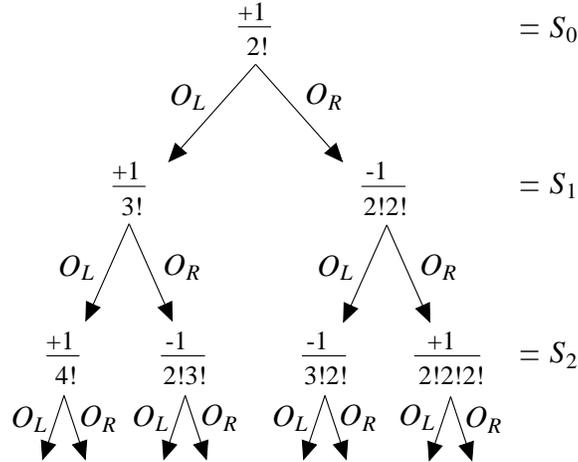}
\caption{The binary tree that generates Bernoulli numbers.}
\end{center}
\end{figure}
The following finite series formed out of the two non-commuting operators
\begin{equation}
S_n = (O_L^{} + O_R^{})^n \!\hf = \left( O_L^n + \sum_{k=0}^{n-1}
O_L^{n-1-k} O_R^{}\,O_L^k + \cdots + O_R^n \right) \!\!\hf .
\label{e:sumtreerep}
\end{equation}
is equivalent to the sum of terms on the $n$-th row of nodes across
the tree.

\begin{theorem}
Bernoulli numbers and the $S_n$ series are related by
\begin{equation}
\fbox{$\; B_n = n!\;S_{n-1} \quad (n\ge 2) \;$}\;.
\label{e:BandS}
\end{equation}
\end{theorem}

For example,
\begin{eqnarray*}
B_3 &=& 3!\, S_2 = 3!\, (O_L^{}+O_R^{})^2 \!\hf 
= 3!\,(O_L^{}+O_R^{})\,(O_L^{}+O_R^{}) \!\hf\\
&=& 3!\,(O_L^{}\,O_L^{} + O_L^{}\,O_R^{} + O_R^{}\,O_L^{} + O_R^{}\,O_R^{}) 
\!\hf\\
&=& 3!\left( \frac{+1}{4!} + \frac{-1}{2!\,3!} + \frac{-1}{3!\,2!} +
\frac{+1}{2!\,2!\,2!} \right) = 0\;.
\end{eqnarray*}

\noindent {\bf Proof}

The Riemann zeta function \cite{titchmarsh}
\begin{equation}
\zeta(s) = \sum_{n=1}^\infty n^{-s} \quad (\re(s)>1, \;s\in\C)
\label{e:zeta}
\end{equation}
can be analytically extended to the left-half of the complex plane
$\re(s)<1$ by the Euler-Maclaurin Summation Formula
(\ref{e:EulerMaclaurinsum}).

Now, instead of adopting (\ref{e:EulerMaclaurinsum}) directly, we
choose to derive the analytic continuation of $\zeta (s)$ into the
left-half of the complex plane step by step.

Consider the difference between (\ref{e:zeta}) with its analogous integral.
\begin{eqnarray*}
\lefteqn{\lim_{N\to\infty} \left[ \sum_{n=1}^N n^{-s} - \left(
\int_1^N x^{-s} dx + \frac{1}{-s\!+\!1} \right) \right]}\\
&=& 1 - \frac{1}{-s\!+\!1} + \lim_{N\to\infty} \sum_{n=2}^N
\left[ n^{-s} - \int_{n-1}^n x^{-s} dx \right]\\
&=& 1 - \frac{1}{-s\!+\!1} + \lim_{N\to\infty} \sum_{n=2}^N 
    \left[ n^{-s} - \frac{1}{-s\!+\!1}\,n^{-s+1} +
    \frac{1}{-s\!+\!1}\,(n-1)^{-s+1} \right]
\end{eqnarray*}

The binomial expansion of the sum over $n$ of the last term, with an
interchange of the order of summation, gives a set of series
$$ 1\,-\,\frac{1}{-s\!+\!1}\,+\,\lim_{N\to\infty}
\left[
\parbox{10cm}{
\begin{eqnarray}
&\bbspace&\bspace\sum_{n=2}^N \left[ n^{-s} - \frac{1}{-s\!+\!1}
 \,n^{-s+1} + \frac{1}{-s\!+\!1}\,n^{-s+1} + (-1)\,n^{-s} \right]
\label{e:series1.1}\\
&\bbspace+&
\frac{(-s)}{2!} \sum_{n=2}^N n^{-s-1} \label{e:series1.2}\\
&\bbspace+&
\frac{(-1)(-s)(-s\!-\!1)}{3!} \sum_{n=2}^N n^{-s-2} \label{e:series1.3}\\
&\bbspace\vdots& \nn\\
&\bbspace+&
\frac{(-1)^{m+1} (-s)(-s\!-\!1)\dots (-s\!+\!1\!-\!m)}{(m\!+\!1)!} 
\sum_{n=2}^N n^{-s-m} \label{e:series1.4}\\
&\bbspace\vdots& \nn
\end{eqnarray}
}
\right]
$$

\begin{eqnarray*}
&=& 1 - \frac{1}{-s\!+\!1} + \lim_{N\to\infty} \sum_{m=1}^\infty 
\left[ \left( \prod_{j=1}^m (-s\!+\!1\!-\!j) \right) \,
\frac{(-1)^{m+1}}{(m\!+\!1)!} \sum_{n=2}^N {-s-m} \right] \\
&=& 1 - \frac{1}{-s\!+\!1} + \lim_{N\to\infty} \sum_{m=1}^\infty 
\left[ \left( \prod_{j=1}^{m} (-s\!+\!1\!-\!j) \right) \,
O_L^{m-1} \!\!\hf \sum_{n=2}^N n^{-s-m} \right]\;,
\end{eqnarray*}
(note that the series (\ref{e:series1.1}) sums to zero). $O_L^{m-1}$
represents a string of $(m-1)$ \, $O_L^{}$ operators, and $O_L^{m-1}
(+1/2!)$ corresponds to the {\em first} term of $S_{m-1}, m\in\Z^+$ in
(\ref{e:sumtreerep}), the {\em first} branch from the left in the
$(m\!-\!1)$-th row of the tree. $O_L^0$ is taken as the identity operator.

Now $\stack{\mbox{\small lim}}{\mbox{\tiny\it N}\to\infty} (1+\int_1^N
kx^{-s} dx = \stack{\mbox{\small lim}}{\mbox{\tiny\it N}\to\infty}
\left[ kN^{-s+1}/(-s\!+\!1) \right]$ converges for $\re(s)>1$, where
arbitrary constants $k,p\in\C$. Therefore, $\stack{\mbox{\small
    lim}}{\mbox{\tiny\it N}\to\infty} \sum_{n=1}^N kn^{-s}$ converges
for $\re(s)>1$. Hence,
$$\zeta(s) = \lim_{N\to\infty}
\left[ \sum_{n=1}^N n^{-s} - \left( \int_1^N x^{-s} dx +
\frac{1}{-s\!+\!1} \right) \right]
= \lim_{N\to\infty} \left[ \sum_{n=1}^N n^{-s}
- \frac{1}{-s\!+\!1} N^{-s+1} \right]$$
defines an analytic continuation of $\zeta(s)$ for $\re(s)>0$,
since the series (\ref{e:series1.2}), (\ref{e:series1.3}) and
(\ref{e:series1.4}) converge for $\re(s)>0$, $\re(s)>-1$, and
$\re(s)>-m\!+\!1$ respectively.

Similarly, find the difference between the series (\ref{e:series1.2})
with its analogous integral.
$$\lim_{N\to\infty} \left[ \sum_{n=2}^N \frac{(-s)}{2!}\,n^{-s-1} 
- \frac{(-s)}{2!} \left( \int_2^N x^{-s-1} dx + \frac{2^{-s}}{-s}
\right)\right]$$
$$=\frac{(-s)\,2^{-s-1}}{2!} \,+\, \frac{2^{-s}}{-s} \,+\, 
\lim_{N\to\infty}
\left[
\parbox{10cm}{
\begin{eqnarray}
&\bbbspace&\bspace
\frac{(-s)(-s\!-\!1)}{2!\,2!}
\sum_{n=3}^N n^{-s-2} \label{e:series2.1}\\
&\bbbspace\vdots&\nn\\
&\bbbspace+& \frac{(-1)^{q+1} (-s)(-s\!-\!1)\dots(-s\!-\!q)}{2!\,(q\!+\!1)!}
\sum_{n=3}^N n^{-s-1-q} \label{e:series2.2}\\
&\bbbspace\vdots&\nn
\end{eqnarray}
}
\right]$$
\begin{eqnarray*}
&\bspace=& \frac{(-s)\,2^{-s-1}}{2!} + \frac{2^{-s}}{-s} 
+ \lim_{N\to\infty} \sum_{m=1}^\infty 
\left[ \left( \prod_{j=1}^{m+1} (-s\!+\!1\!-j) \right)
\,\frac{(-1)^{m+1}}{2!\,(m\!+\!1)!} \sum_{n=3}^N n^{-s-m-1} \right]\\
&\bspace=& \frac{(-s)\,2^{-s-1}}{2!} + \frac{2^{-s}}{-s} 
+ \lim_{N\to\infty} \sum_{m=1}^\infty 
\left[ \left( \prod_{j=1}^{m+1} (-s\!+\!1\!-\!j) \right) \,
O_R^{}\,O_L^{m-1} \!\!\hf \sum_{n=3}^N n^{-s-m-1} \right]
\end{eqnarray*}
where $O_R^{}\,O_L^{m-1} (+1/2!)$ corresponds to the {\em second} term of
$S_m, \;m\in\Z^+$, the {\em second} branch from the left in the $m$-th
row of the tree, and $q\in\Z^+, \;q>2$.  Hence,
\begin{eqnarray*}
\zeta(s)&=& \lim_{N\to\infty}
\left[ \sum_{n=1}^N n^{-s} - \left( \int_1^N x^{-s} + \frac{1}{-s+1} dx
\right)
- \frac{(-s)}{2!} \left( \int_2^N x^{-s-1} + \frac{2^{-s}}{-s} dx
\right)\right] \\
&=& \lim_{N\to\infty} \left[ \sum_{n=1}^N n^{-s}
- \frac{1}{-s\!+\!1}\,N^{-s+1} - \frac{1}{2!}\,N^{-s} \right] \\
&=& \lim_{N\to\infty} \left[ \sum_{n=1}^N n^{-s}
- \frac{1}{-s\!+\!1}\,N^{-s+1} - S_0\,N^{-s} \right]
\end{eqnarray*}
defines an analytic continuation of $\zeta(s)$ for $\re(s)>-1$, since
the series (\ref{e:series2.1}) and (\ref{e:series2.2}) converge for
$\re(s)>-1$ and $\re(s)>-q$ respectively. Note that $(+1/2!)$, the
coefficient of $N^{-s}$, corresponds to $S_0$, the starting node
(zeroth row) of the tree.

To analytically extend $\zeta(s)$ to $\re(s)>-2$, we again subtract
in a similar way the set of all the series of $n^{-s-2}$ (e.g.,
(\ref{e:series1.2}) and (\ref{e:series1.3})), which diverge for
$\re(s)\leq -1$, to get sets of series of higher orders that
converge for $\re(s)>-2$.

Similarly, the difference between the series (\ref{e:series1.2}) with
its analogous integral, and the corresponding one for series
(\ref{e:series1.3}) are
$$\frac{(-s)(-s\!-\!1)\,3^{-s-1}}{2!\,2!} \,+\, \frac{3^{-s-1}}{-s-1}
\,+\,\lim_{N\to\infty} \sum_{m=1}^\infty 
\left[ \left( \prod_{j=1}^{m+3} (-s\!+\!1\!-\!j) \right) \,
O_L^{}\,O_R^{}\,O_L^{m-1} \!\!\hf \sum_{n=4}^N n^{-s-m-3} \right]$$
and
$$\frac{(-s)(-s\!-\!1)\,3^{-s-1}}{2!\,2!} \,+\, \frac{3^{-s-1}}{-s-1}
\,+\,\lim_{N\to\infty} \sum_{m=1}^\infty 
\left[ \left( \prod_{j=1}^{m+3} (-s\!+\!1\!-\!j) \right) \,
O_R^2\,O_L^{m-1} \!\!\hf \sum_{n=4}^N n^{-s-m-3} \right]$$
respectively, where $O_L^{}\,O_R^{}\,O_L^{m-1} (+1/2!)$ and
$O_R^2\,O_L^{m-1} (+1/2!)$ correspond to the {\em third} and {\em
  fourth} terms respectively of $S_{m+1}, m\in\Z^+$, the {\em third}
and {\em fourth} branch from the left in the $(m\!+\!1)$-th row of the
tree. Hence,
\begin{eqnarray}
\zeta(s) &=& \lim_{N\to\infty}
\left[ \begin{array}{l}
\ds \sum_{n=1}^N n^{-s} - \left( \int_1^N x^{-s} dx + \frac{1}{-s+1} 
\right)\\
\ds\Bigg. {}- \frac{(-s)}{2!} \left( \int_2^N x^{-s-1} dx +
\frac{2^{-s}}{-s} \right)\\
\ds\Bigg. {}- (-s)(-s\!-\!1)\,O_L^{} \!\hf \left(\int_3^N x^{-s-2} dx 
+ \frac{3^{-s-1}}{-s-1} \right)\\
\ds\Bigg. {}- (-s)(-s\!-\!1)\,O_R^{} \!\hf \left(\int_3^N x^{-s-2} dx
+ \frac{3^{-s-1}}{-s-1} \right) \end{array} \right]\nn\\
&=& \lim_{N\to\infty} 
\left[ \begin{array}{l}
\ds \sum_{n=1}^N n^{-s} - \frac{1}{-s\!+\!1}\,N^{-s+1} - S_0\,N^{-s} 
- S_1\,(-s)\,N^{-s-1}
\end{array}
\!\!\!\right]\label{e:seriesorder2}
\end{eqnarray}
defines an analytic continuation of $\zeta(s)$ for $\re(s)>-2$, where
$S_0 = (+1/2!)$ and $S_1 = (O_L^{} + O_R^{})\,(+1/2!)$.

Comparing series (\ref{e:seriesorder2}) with the Euler-Maclaurin
Summation Formula (\ref{e:EulerMaclaurinsum}) for $M=1$,
$$\zeta(s) = \lim_{N\to\infty} \left[ \sum_{n=1}^N n^{-s}
- \frac{1}{-s\!+\!1}\,N^{-s+1} - \frac{1}{2}\,N^{-s}
- \frac{B_2}{2!}\,(-s)\,N^{-s-1} \right]\;,$$
we find that $B_2 = 2!\,S_1$.

At this point, we observe that every time we find the difference
between a divergent series and its analogous integral, the resulting
binomial expansion has the effect of inserting additional operators
$O_R^{}\,O_L^{m-1}$ immediately before $(+1/2!)$. When we write down these
sequences, the tree appears.

To obtain further analytic continuation of $\zeta(s)$, we make use of
the tree and write
\begin{eqnarray*}
\zeta(s) &=& \lim_{N\to\infty}
\left[ \begin{array}{l}
\ds \sum_{n=1}^N n^{-s} - \frac{1}{-s\!+\!1}\,N^{-s+1} - S_0\,N^{-s} \\
\ds\Bigg. -\,\sum_{m=1}^M \left[ \left( \prod_{j=1}^{m} (-s\!+\!1\!-\!j)
  \right) S_m\,N^{-s-m} \right]
\end{array} \right]
\end{eqnarray*}
where $\ds\qquad S_m = (O_L^{} + O_R^{})^m \!\hf$.

Comparing the above with the Euler-Maclaurin Summation Formula
(\ref{e:EulerMaclaurinsum}), and noting that $B_n$ vanishes for odd
$n\geq 3$, and so the coefficients of $N^{-s-n+1}$ also vanish for odd
$n \geq 3$, we get (\ref{e:BandS}).

\hfill\qed

By observation, the sum-across-the-tree representation of $S_n$ in
(\ref{e:sumtreerep}) can also be seen to be equivalent to the following
determinant known to generate $B_n$.
$$
S_n = (-1)^n \left| \begin{array}{ccccccc}
\ds\Bigg.\frac{1}{2!}&1&0&\;\quad 0\quad\;&\;\quad 0\quad\;&
 \;\cdots\;&\;\;\;0\;\\
\ds\Bigg.\frac{1}{3!}&\ds\frac{1}{2!}&1&0&0&\cdots&\;\;\;0\;\\
\ds\Bigg.\frac{1}{4!}&\ds\frac{1}{3!}&\ds\frac{1}{2!}&1&0&\cdots&\;\;\;0\;\\
\ds\Bigg.\vdots&\ddots&\ddots&\ddots&\ddots&\ddots&\;\;\;\vdots\;\\
\ds\Bigg.\frac{1}{(n-2)!}&\ddots&\ddots&\ds\frac{1}{3!}&\ds\frac{1}{2!}&1&
 \;\;\;0\;\\
\ds\Bigg.\frac{1}{(n-1)!}&\ds\frac{1}{(n-2)!}&\ddots&\ddots&
 \ds\frac{1}{3!}&\ds\frac{1}{2!}&\;\;\;1\;\\
\ds\Bigg.\frac{1}{n!}&\ds\frac{1}{(n-1)!}&\ds\frac{1}{(n-2)!}&
 \ddots&\ddots&\ds\frac{1}{3!}&\;\;\;\ds\frac{1}{2!}\; \end{array}
\right| .
$$

\section{The Tree-Generating Operator and Bernoulli Function}

We can now expand the tree-generating operator $(O_L + O_R)$ raised to the
complex power $(s\!-\!1)$ acting on $\ds\hf$ as follows.

\begin{lemma}
\label{l:optreeseries}
\begin{eqnarray}
\lefteqn{}\label{e:optreeseries}\\
&&\bspace (O_L^{} + O_R^{})^{s-1} \! \hf\nn\\
&=& \left( \Big. w \id - \left[ \big. w \id - (O_L^{} + O_R^{}) \right]
   \right)^{\!s-1} \hf\nn\\
&=& w^{s-1}\Big(\id - \Big[\id - \frac{1}{w} (O_L^{} + O_R^{}) 
   \Big]\Big)^{\!s-1} \hf\nn\\ 
&=& w^{s-1}\!\left( \id + \sum_{n=1}^\infty (-1)^n
   {s-1 \choose n}\!\left[ \big. \id - \frac{1}{w} (O_L^{} + O_R^{}) 
   \right]^{\!n} \right) \!\hf\nn\\ 
&=& w^{s-1}\left(\id + \sum_{n=1}^\infty (-1)^n {s\!-\!1 \choose n} \!
\left[ \id + \sum_{m=1}^n 
    \left(\frac{-1}{w}\right)^{\!\!m} \!\! {n \choose m}
    (O_L^{} + O_R^{})^m \right] \right) \!\hf \nn\\
&=& w^{s-1}\left(\frac{1}{2} + \sum_{n=1}^\infty 
    (-1)^n {s\!-\!1 \choose n} \! \left[ \frac{1}{2} +
    \sum_{m=1}^n \left(\frac{-1}{w}\right)^{\!\!m} \!\! {n \choose m}
    \frac{B_{m+1}}{(m\!+\!1)!}\right]\right)\nn
\end{eqnarray}
which converges for $\re(s)>(1/w)$ where $s\in\C, \;w\in\R, \;w>0$,
and $\id$ is the identity operator.
\end{lemma}

\begin{theorem}[{\bf Bernoulli Function}]
\begin{eqnarray}
B(s) &\!=\!& \Gamma(1\!+\!s) \; (O_L^{} + O_R^{})^{s-1} \! \hf
 \label{e:Bseries}\\
&\!=\!& w^{s-1}\;\Gamma(1\!+\!s) \left(
 \frac{1}{2} + \sum_{n=1}^\infty (-1)^n {s\!-\!1 \choose n} \! \left[
 \frac{1}{2} + \sum_{m=1}^n \left(\frac{-1}{w}\right)^{\!\!m} \!\! 
{n \choose m} \frac{B_{m+1}}{(m\!+\!1)!} \right] \right)\nn
\end{eqnarray}
which converges for $\re(s)>(1/w)$ where $s\in\C, \;w\in\R, \;w>0$.
\end{theorem}

\noindent{\bf Proof}

From (\ref{e:sumtreerep}) and (\ref{e:BandS}), we have
$$B_n = n! \; (O_L^{} + O_R^{})^{n-1} \!\hf = \Gamma(1\!+\!n) \,
(O_L^{} + O_R^{})^{n-1} \!\hf \quad(n\ge 2)\;.$$

Analytically extending the tree-generating operator $(O_L^{} + O_R^{})$
with (\ref{e:optreeseries}) in Lemma \ref{l:optreeseries} effectively
turns the sequence of $B_n$ into a function $B(s)$ as the analytic
continuation of $B_n$.

\hfill\qed

\begin{figure}[hbt]
\begin{center}
\includegraphics[height=157pt,width=250pt]{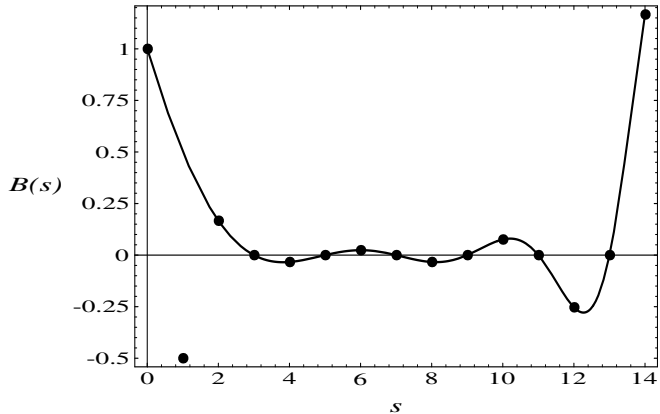}
\caption{The curve $B(s)$ runs through the points of all $(n,B_n)$
  except $(1,B_1)$.\label{fig:B(s)}}
\end{center}
\end{figure}

Figure \ref{fig:B(s)} shows a plot of $B(s)$ for real $s$. All the
Bernoulli numbers $B_n$ agree with Bernoulli function $B(n)$ except at
$n=1$, i.e., $B(n)=B_n$ for $n=0$ or $n\ge 2$,
\begin{equation}
B(1) = 1/2 \quad \mbox{ but } \quad B_1 = -\,1/2\;.
\end{equation}
We shall now address the surprising discrepancy between $B(1)$ and $B_1$.

\section{Fixing of the Arbitrary Sign Convention of $B_1$}

$B_1=-1/2$ has largely been adopted as the standard sign convention
partly due to elegance in notation and partly due to its widespread
usage although there had been suggestions for favoring the sign
convention $B_1=1/2$.

Looking back at (\ref{e:defBn}) to (\ref{e:defBnrecursion}) in Definition
\ref{d:Bn}, we see that the sign convention of $B_1$ was arbitrary.
Figure \ref{fig:B(s)} show that the analytic continuation of $B_n$
actually fixes the arbitrary sign convention of $B_1$. A mathematical
fact should precede notational elegance or personal preference.

\begin{definition}[{\bf Redefinition of Bernoulli Numbers}]\quad\\
\label{d:redefB}
To have consistency between Bernoulli numbers $B_n$ and their analytic
continuation $B(s)$, we should redefine $B_n$ as
\begin{equation}
\frac{z}{e^z-1} = \sum_{n=0}^\infty (-1)^n \frac{B_n}{n!} z^n
\quad (|z|<2\pi)\;,\label{e:newdefBn} \end{equation}
or
\begin{equation}
B_n = \frac{(-1)^{n+1}}{n\!+\!1} \,\sum_{k=0}^{n-1}\,(-1)^k\, 
{n\!+\!1 \choose k}\, B_k\;, \quad B_0 = 1\;.
\label{e:newdefBnrecursion} \end{equation}
\end{definition}
The factor $(-1)^n$ introduced in (\ref{e:newdefBn}) and
(\ref{e:newdefBnrecursion}) only changes the sign in the conventional
definition of the only non-zero odd Bernoulli numbers, $B_1$, from
$B_1 = -1/2$ to the redefined $B_1 = B(1) = 1/2$.

\section{Proof of the New Representation of $\zeta(s)$}

We now have covered sufficient concepts to prove Theorem \ref{t:zetarep}.

{\bf Proof}

$B_n$ are related to the Riemann zeta function $\zeta(s)$ as
\cite[p.34]{bateman}
\begin{eqnarray*}
\zeta(-n)&=& - \;\frac{B_{n+1}}{n\!+\!1} \qquad\qquad\qquad (n\in\Z^+)\;,\\
\zeta(2n)&=& \frac{(-1)^{n+1}(2\pi)^{2n}}{2\;(2n)!} \,B_{2n} 
\quad (n\in\Z^+)\;.
\end{eqnarray*}

Hence, $\zeta(1\!-\!s) = - \;B(s)/s \;\;(s\in\C)$.  Replacing $B(s)$
with the series in (\ref{e:Bseries}) and noting that
$\Gamma(1\!+\!s)/s = \Gamma(s)\;$ gives
\begin{equation}\zeta(1\!-\!s) = -\;w^{s-1}\;\Gamma(s) \left( 
  \frac{1}{2} + \sum_{n=1}^\infty (-1)^n  {s\!-\!1 \choose n} \!
  \left[ \frac{1}{2} + \sum_{m=1}^n \left(\frac{-1}{w}\right)^{\!\!m} 
  \!\! {n \choose m} \frac{B_{m+1}}{(m\!+\!1)!} \right] \right)
\label{e:zeta(1-s)series}
\end{equation}
which converges for $\re(1\!-\!s) < 1-(1/w)$ where $s\in\C, \;w\in\R,
 \;w>0$.

The functional equation of the Riemann zeta function \cite{ahlfors}
relates $\zeta(1\!-\!s)$ to $\zeta(s)$ as
\begin{equation}
\zeta(1\!-\!s) = 2 \;(2\pi)^{-s} \;\Gamma(s)\; \cos\left(\frac{\pi
    s}{2}\right) \;\zeta(s)\;.
\label{e:funceqn}
\end{equation}

Applying the functional equation (\ref{e:funceqn}) to
(\ref{e:zeta(1-s)series}) yields
\begin{equation}\cos\left(\frac{\pi \hat{s}}{2}\right)
\zeta(s) = -\; \frac{(2\pi)^s}{2} (w)^{s-1} \left( \frac{1}{2} +
  \sum_{n=1}^\infty (-1)^n {s\!-\!1 \choose n}
    \!  \left[ \frac{1}{2} + \sum_{m=1}^n
    \left(\frac{-1}{w}\right)^{\!\!m} \!\!  {n \choose m}
    \frac{B_{m+1}}{(m\!+\!1)!} \right] \right)
\label{e:zetacos}
\end{equation}
(\ref{e:zetacos}) in the limit form gives the Theorem.

\hfill\qed

\end{document}